# Automation of Trimming Die Design Inspection by Zigzag Process Between AI and CAD Domains

Jinsub Lee[†]

Department of Artificial Intelligence, Sungkyunkwan University, Suwon, Gyeonggi-do, 16419, Gyeonggi-do, Republic of Korea

Tae-Hyun Kim[†]

Department of Smart Fab. Technology, Sungkyunkwan University, Suwon, Gyeonggi-do, 16419, Gyeonggi-do, Republic of Korea

Sang-Hwan Jeon

Stamping Tool Design Team, Kia motors, Hwaseong, Gyeonggi-do, 18571, Republic of Korea

Sung-Hyun Park

Stamping Tool Design Team, Kia motors, Hwaseong, Gyeonggi-do, 18571, Republic of Korea

Sang-Hi Kim

Advanced Manufacturing CAE Team, Hyundai motors, Hwaseong, Gyeonggi-do, 18280, Republic of Korea

Eun-Ho Lee[*]

Department of Smart Fab. Technology, Sungkyunkwan University, Suwon, Gyeonggi-do, 16419, Gyeonggi-do, Republic of Korea; and

School of Mechanical Engineering, Sungkyunkwan University, Suwon, Gyeonggi-do, 16419, Republic of Korea

Jee-Hyong Lee[*]

Department of Artificial Intelligence, Sungkyunkwan University, Suwon, Gyeonggi-do, 16419, Gyeonggi-do, Republic of Korea

† Equal contributions

* Corresponding author. Tel.: +82-31-290-7987; fax: +82-31-299-4637.(J.-H. Lee), Tel.: +82-31-299-4861.(E.H. Lee)

*E-mail address:* john@skku.edu(J.-H. Lee), e.h.lee@skku.edu(E.H. Lee)

[Lee et al., 2022, Automation of Trimming Die Design Inspection by Zigzag Process Between AI and CAD Domains]

**Abstract**

Quality control using artificial intelligence (AI) has led to considerable improvements in the manufacturing industry. However, trimming die design inspection to minimize potential failures is still done manually by engineers. The inspection is labor-intensive, time-consuming, and prone to human error. In this study, we developed an automatic design inspection system for automobile trimming dies by integrating AI modules and computer-aided design (CAD) software. The AI modules replace the engineers' judgment, and the CAD software performs the operations requested by the AI modules. The design inspection is completed through a zigzag interaction between the AI modules and CAD software by a one-click operation without the intervention of experts. In the zigzag interaction process, 3D CAD files are converted to 2D images for interaction from the CAD software to the AI modules, and the output of the AI modules is transferred to the CAD software through VBA macros. The AI modules are designed to be CAD-independent and data-efficient; therefore, they can be easily adapted to other CAD software and achieve a high performance, even when trained with only a few trimming dies. The length measurement error is only 2.4 % on an average, which is at the 14th place out of 18 grades according to ISO 286-1:1988 standard tolerance. In addition, the inspection time is reduced to approximately one-fifth of the manual inspection time by experts.

[Lee et al., 2022, Automation of Trimming Die Design Inspection by Zigzag Process Between AI and CAD Domains]

# 1. Introduction

In recent years, manufacturing quality control has seen considerable improvements. Throughout the manufacturing process, product reviews and feedback are essential for quality improvement. Among the various methods for quality improvement, failure mode and effects analysis (FMEA) has been widely used in manufacturing. FMEA is used for the identification, tracking, and impact analysis of potential failure modes to reduce the defects in processes and products (Sader et al., 2020). Engineers minimize potential failures by analyzing previous failure modes (Wu et al., 2021). However, because most FMEA processes rely on engineers, they are labor-intensive, time-consuming, and prone to human errors.

Several studies have presented automated FMEA by utilizing large amounts of data and artificial intelligence (AI), which can handle the high dimensionality, complexity, and dynamics of the manufacturing industry. Chang et al. proposed an evolutionary tree learning algorithm-based approach for clustering and visualizing the failure modes in FMEA (Chang et al., 2015). Jensen et al and Boral et al. identified potential failures by using AI (Jensen et al., 2014; Boral et al., 2019). Geren et al. (2017) proposed a new 3D design platform to automate the repetitive parametric part design and assembly process for various types of ball-joint assemblies. Long et al. (2021) proposed a system for automating the design of mechanical products without expert participation. FMEA automation using AI is performed in various processes; however, automatic design inspection has rarely been studied.

Automation of product design inspection, which is a feasibility review to determine the validity of design prototypes, is beneficial to a manufacturing process because it can significantly reduce the waste of resources and costs by effectively detecting early failure modes. In a design process, the design modeling time is significantly reduced by an automatic design based on parametric 3D modeling. However, the design inspection process is still manually performed by engineers and requires considerable resources and time.

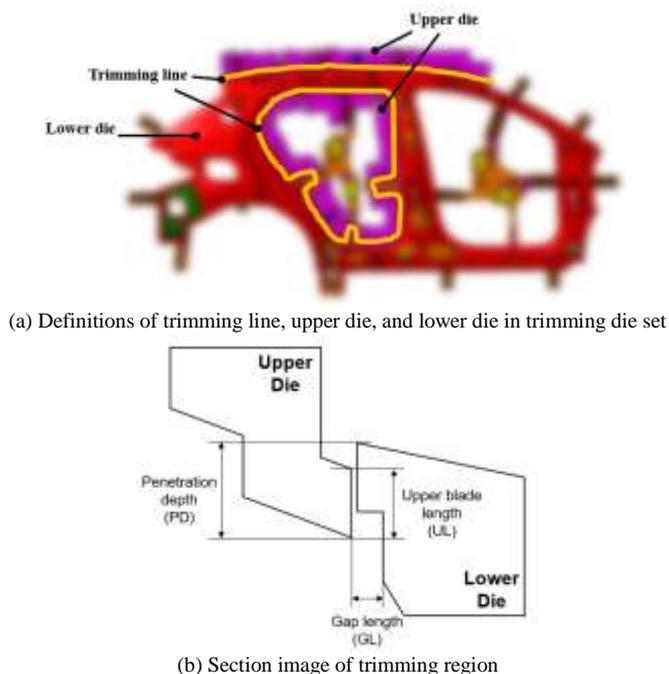

(a) Definitions of trimming line, upper die, and lower die in trimming die set

(b) Section image of trimming region

Fig. 1. Inspection target. (a) The trimming lines are orange lines which should cut the excess parts of the steel sheet after forming processes. (b) PD, UL, and GL are defined. Trimming sections should satisfy PD, UL, and GL to cut designated parts clearly with balanced across the trimming line. Fig. 1(a) is partially blurred by the manufacturer's requests.

This study aims to develop a system that can automatically inspect trimming die designs using AI. Fig. 1(a) shows the design of an automobile trimming die. A trimming die set is a mechanical tool used to trim excess parts of a product body by punching the lower die with the upper die, and a trimming line is used to cut the spare part. Note that the CAD drawing of the die is blurred throughout the manuscript, as shown in Fig. 1(a), with respect to the manufacturer's request. An illustration of a cross-section perpendicular to the trimming line is shown in Fig. 1(b). The main design variables of the trimming section include penetration depth (PD), upper blade length (UL), and gap length (GL). PD is the displacement of the upper die and UL is the length of the blade of the upper die. These are important parameters for the trimming process to be controlled as designed. GL is the interval between the blade and the body of the lower die and is required to facilitate the removal of excess materials without stock.

[Lee et al., 2022, Automation of Trimming Die Design Inspection by Zigzag Process Between AI and CAD Domains]

The inspection of trimming dies is important because it is related to productivity and safety. Trimming dies need to be designed such that scraps are properly extracted during the trimming process. If the trimming dies are improperly designed, the scraps may not be extracted, leading to faulty (defective) dies. However, manual design inspection of trimming dies is a labor-intensive task, which reduces productivity. Engineers must repeatedly check whether the design meets certain standards by considering sections at every inspection spot.

We develop an automatic design inspection system for trimming dies by integrating AI modules into a computer-aided design (CAD) software. The AI modules and CAD software are integrated into a zigzag interaction process. The AI modules replace the engineers' judgment, and the CAD software performs the operations requested by the AI modules. The design inspection is completed through automatic feedback between the AI modules and the CAD software in the zigzag interaction process. Because 3D CAD files are converted into 2D images in the zigzag process, the AI modules can run with any CAD software; that is, they can be easily adapted to other CAD software. The AI modules are data-efficient and designed to achieve a high performance, even when trained with only a few trimming dies.

Our inspection system consists of data processing, trimming region detection, selective cropping, and length measurement stages. In the data processing stage, sectioned image files are extracted from the CAD assembly files by excluding the background and unnecessary geometry. The trimming regions are identified in the sectioned images using a deep object detection model in the trimming region detection stage. The enlarged images of the identified trimming regions are then extracted together with scale factors in the selective cropping stage. In the target point detection and length measurement stage, the target line endpoints are detected, and the lengths are measured using scale factors.

Our system measures the dimensions of PD, UL, and GL from 3D CAD assembly files using 2D convolutional neural network (CNN) models. The measuring process is fully automated by a one-click operation without expert intervention. It can reduce the inspection time and cost with accurate length measurements. The inspection time is reduced to approximately one-fifth of the manual inspection time by experts. The length measurement error is only 2.4% on average. In addition, this CNN-based method allows the developed AI module to be connected regardless of the type of CAD program.

In the following section, related studies using AI in the manufacturing process and the object detection in our automatic inspection system are introduced. In Section 3, the proposed CAD-based automatic inspection system is described. In Sections 4 and 5, we present the experimental results, analysis, and discussion. Finally, we conclude the paper in Section 6.

## 2. Related Work

Automation using AI has been widely adopted in various manufacturing fields. We discuss this from the viewpoint of intelligent manufacturing and automation. We also introduce Faster R-CNN (Ren et al., 2015), which is employed in our system to identify trimming regions and target line endpoints.

*2.1. Intelligent manufacturing technology with AI*

**Intelligent manufacturing process** is closely related to applications in smart factories, which adopt a combination of physical and cyber technologies (Brettel et al., 2014; Li et al., 2016a; Chen et al., 2017; Zhou et al., 2018). Intelligent manufacturing technology has emerged and is reported to have a profound and lasting effect on advanced manufacturing technology using the acquired data (Lucke et al., 2008). However, different formats of raw data from sensors in the production line has created technical problems in the gap between manufacturing and information fields, which include high dimensionality, complexity, and dynamicity (Loyer et al., 2016; Davis et al., 2015; Wuest el al., 2016). Many studies have reported the adoption of artificial neural networks (ANNs) to overcome high complexity, which involves both physics and practical factors in the real implementation of assessments in manufacturing industries (Boral et al., 2019; Li et al., 2013). To optimize a manufacturing system consisting of a machine, robot, conveyor, maintenance support, and material handling system, Li et al. (2016b) proposed a production process control of the system throughput for a multiproduct production line by considering a rework loop with AI. Kumar (2017) reviewed various cases where AI was used in computer-aided process planning and manufacturing. Early diagnosis and monitoring of product quality through quality control has been applied to AI (Wang et al., 2019). The detection of surface defects on the product in production lines also uses AI to ensure accuracy (Yang et al., 2005), along with maintenance scheduling for machines and robots in manufacturing lines (Kumar et al., 2018). The problem of dynamic dispatching for unreliable machines in re-entrant semiconductor production systems was considered in the study by Wu et al. (2020). They combined a deep neural network model and Markov decision processes (MDPs) to rapidly generate near-optimal dynamic control policies for problems that are too large to be solved only by MDPs, thus showing the potential of machine learning in controlling unreliable manufacturing systems (Sidford et al., 2018; Chien et al., 2020). Elhoone et al. (2020) proposed three ANN algorithms to support the dynamic allocation of digital designs to different additive manufacturing techniques.

**Process automation** is crucial for implementing mass production in the manufacturing industry through a meta-database of pre-processed sensor data. As mentioned in the Introduction, Geren et al. (2017) proposed a new platform that allowed the automated reconfigurations of a ball-joint assembly model. Long et al. (2021) proposed a knowledge-based automated design system to minimize the in-person participation and interaction of experts. A design inspection process is required to check whether several requirements are satisfied, which is based on FMEA. Automated visual inspection processes using high-resolution cameras, such as systems for automated power line inspection (Martinez et al., 2018) and bean quality inspection (De Araújo et al., 2015) using AI, are currently being used to capture defects and cracks. Similarly, an inspection must be performed within the CAD interface to check if design requirements are satisfied in cyberspace, as presented in this work. AI can interact with CAD systems (Shah et al., 2001; Engelmore et al., 1993; La Rocca, 2012; Pokojski et al., 2022), automate repetitive and non-creative design tasks, and support a multidisciplinary design optimization in all phases of the design process (Engelmore et al., 1993). AI can be applied through an inference mechanism that can access knowledge storage structured and formalized using knowledge representation, such as rules and frames (Biedermann et al., 2021). FMEA knowledge can also be stored in a structured space to facilitate the work of experts. Therefore, by using knowledge-based engineering with AI, it is possible to develop an automatic inspection application that can imitate human experts and automate the problem-solving process, as presented in this study.

*2.2. Faster R-CNN: an object detection method*

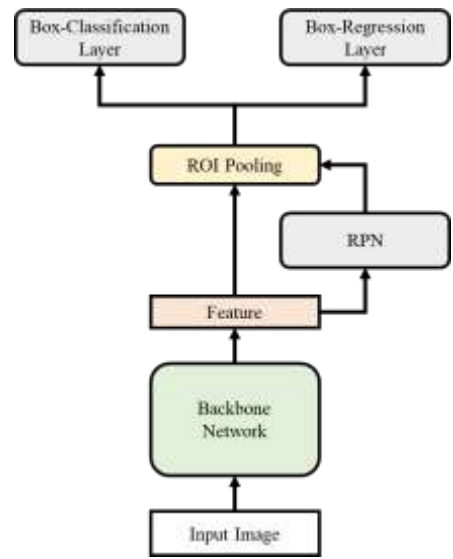

Fig. 2. Diagram of Faster R-CNN. RPN is the region proposal network. Faster R-CNN inherits the Fast R-CNN structure, improves both speed and accuracy by removing selective search and adding RPN.

Object detection is an important task in computer vision to find regions where objects are likely to exist in an image (region proposal) and identify the object classes (classification). These can be classified into one- and two-stage detections. YOLO (Redmon et al., 2016), SSD (Liu et al., 2016), and Retina-Net (Lin et al., 2017) are categorized as one-stage detectors, whereas two-stage detectors include R-CNN (Girshick et al., 2014), Fast R-CNN (Girshick, 2015), and Faster R-CNN (Ren et al., 2015). One-stage detectors, which simultaneously perform region proposals and classifications, are relatively faster, but have a low accuracy. In contrast, two-stage detectors, which perform sequential region proposals and classifications, are relatively slow, but more accurate. This study uses a two-stage detector because accuracy is important for item assessment in the manufacturing industry. Nevertheless, an automated inspection using two-stage detectors is still considerably faster than human inspection.

The proposed automated inspection system uses Faster R-CNN, which inherits the Fast R-CNN structure to accurately detect trimming sections and measure the target lengths. Faster R-CNN is a method to increase both speed and accuracy by generating a region of interest (ROI) through a region proposal network (RPN). As shown in Fig. 2, Faster R-CNN includes a backbone network, RPN, ROI pooling, box regression layer, and box-classification layer. The backbone network extracts feature maps from the input images and the RPN generates region proposals from the feature maps. The training loss of RPN in Faster R-CNN is given as:

$$l_{rpn}(\{p_i\},\{t_i\}) = \frac{1}{N_{cls}}\sum_i l_{cls}(p_i, p_i^*) + \lambda \frac{1}{N_{reg}}\sum_i p_i^* l_{reg}(t_i, t_i^*) \quad (1)$$

The RPN in Faster R-CNN is trained with a classification loss over two classes (object vs. non-object) and a bounding box regression loss (smooth L1 loss) (Girshick, 2015). In Eq. (1), $i$ is the index of an anchor, which is a box used as a candidate for objects in a mini-batch, and $p_i$ is the predicted probability that anchor $i$ is an object. We denote the ground-truth label as $p_i^*$, and the predicted bounding box and ground-truth box as $t_i$ and $t_i^*$, respectively.

ROI pooling is applied to the generated fixed-size feature maps using region proposals and feature maps. The classification and box regression are then trained using fixed-size feature maps. Eq. (2) gives the multitask loss of Fast R-CNN. A cross-entropy loss is used for classification and smooth L1 loss for box regression.

$$l(p, u, t^u, v) = l_{cls}(p, c) + \lambda[u \geq 1]l_{loc}(t^u, v) \quad (2)$$

Ground-truth data is labeled with class $c$ and box regression target $v$. The discrete probability distribution over k categories is $p = (p_1, \ldots, p_k)$. The probability distribution $p$ is computed by a softmax function over the k outputs of the classification layer, and the box regression offsets for each of the K object classes are $t^k = (t_x^k, t_y^k, t_w^k, t_h^k)$. The log loss for a true class $c$ is $l_{cls}(p, c) = -log\ p_c$. For box regression, the loss function is

$$l_{loc}(t^u, v) = \sum_{i \in \{x,y,w,h\}} \text{smooth}_{L_1}(t_i^u - v_i), \quad (3)$$

where

$$\text{smooth}_{L_1}(x) = \begin{cases} 0.5x^2 & \text{if } |x| < 1 \\ |x| - 0.5 & \text{otherwise.} \end{cases} \quad (4)$$

## 3. Proposed Method: Automatic Design Inspection

*3.1. Overview*

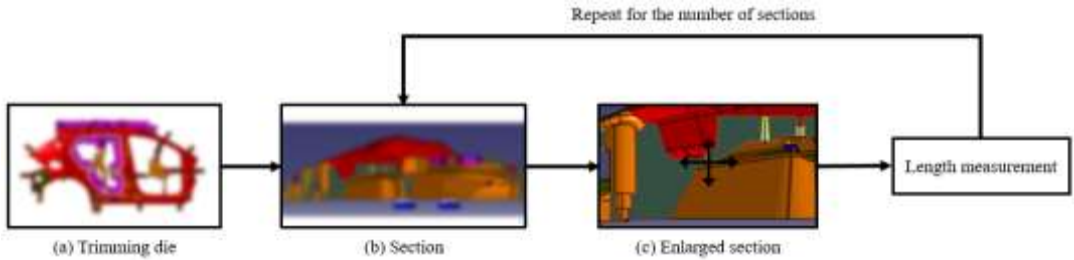

Fig. 3. Manual trimming die inspection process. Fig. 3(b) is partially blurred by the manufacturer's requests.

Fig. 3 shows the process of trimming die inspection by engineers. The trimming die inspection task is labor-intensive. First, engineers choose a trimming line, as shown in Fig. 3(a), and generate a section image perpendicular to it using the sectioning function built in CAD. The section images are shown in Fig. 3(b). The engineers then zoom into the trimming region for an inspection, as shown in Fig. 3(c), and measure the lengths of the target lines, such as PD, UL, and GL. Finally, they determine the validity of the measured lengths based on FMEA guidelines. This process is repeated for all sections to be inspected, which is labor-intensive. However, it is difficult to automate the inspection of trimming dies because of the numerous interactions between the engineers and CAD software.

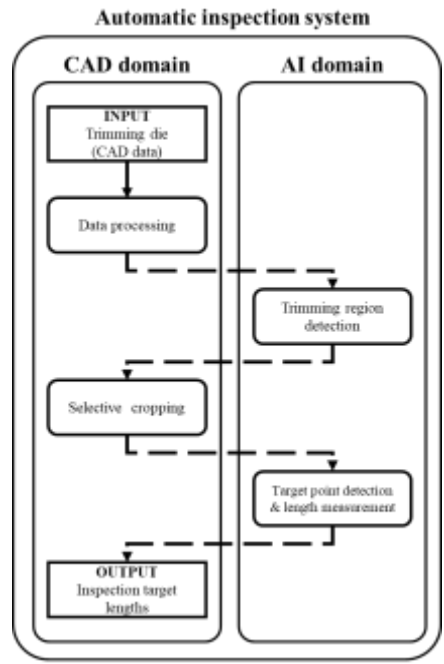

Fig. 4. Structure of the proposed automatic design inspection system. Dashed lines represent the zigzag interaction between CAD and AI domains. The inspection process is composed of a series of zigzag interaction between CAD and AI domains.

To overcome this difficulty, we propose a zigzag interaction process between the AI and CAD domains. The flow of the automated inspection is illustrated in Fig. 4. Our system consists of CAD domain operations—data processing and selective cropping, and AI domain operations—trimming region detection, target point detection, and length measurement. The AI modules generate decisions on behalf of the engineers, which are then fed back to the CAD domain. Then, the modules in the CAD domain use them to generate the data required by the AI modules for the next step.

Data processing in the CAD domain, which is implemented using the macro function of CAD, creates 2D trimmed section images from 3D CAD files for the AI models. In this process, CAD data processing, which includes background removal, color reduction, and tree information removal, is performed, and shortcut marks are added such that the 2D CNN can easily learn the trimming region where the target line is located. Then, a trimming region detection in the AI domain detects trimming regions in the section images using the 2D CNN. Selective cropping in the CAD domain extracts enlarged trimming region images from the CAD file based on the results of the trimming region detector. Target point detection and length measurement in the AI domain detects both endpoints of the target lines in the enlarged trimming region images using Faster R-CNN and measures the target lengths. To smoothly automate the operations in CAD and AI domains, we tightly integrate the CAD software with AI models.

Four trimming die CAD files, supplied by the Hyundai-KIA automotive group, are used to develop the automatic inspection system. The data processing generates a total of 200 section images, including 50 section images from each CAD file. The image size is 1,880 × 933 pixels.

*3.2. Data processing stage*

In data processing, section images are generated as inputs to the trimming region detector, as shown in Fig. 5. The data processor design is based on the CAD's macro functions. First, the engineers chose a trimming line for inspection and determine the inspection interval. The inspection spots are marked at every inspection interval on the chosen trimming line, as shown in Fig. 6(a). The data processor then automatically generates section images along the sectioning line to the trimming line at each inspection spot, as shown in Fig. 6(b). We align the section images into a front view using the normal vector of the section, as shown in Fig. 6(c). We then perform the background removal, color reduction, and tree information removal using the CAD's built-in function to generate simplified images, as shown in Fig. 6(d). We set the CAD background to white and lines to black and fill the polygons with gray color. Section images at every inspection spot are automatically generated by the CAD's built-in capturing tool.

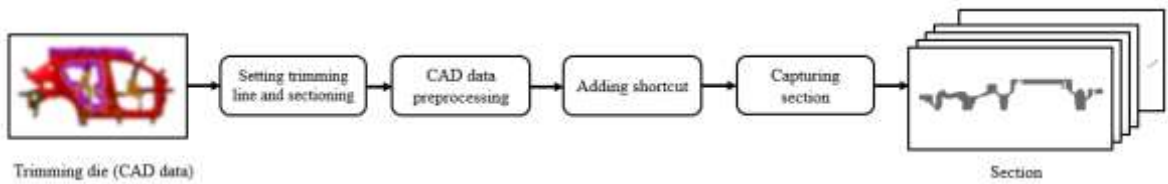

Fig. 5. The process of data processing. The trimming die image is partially blurred by the manufacturer's requests.

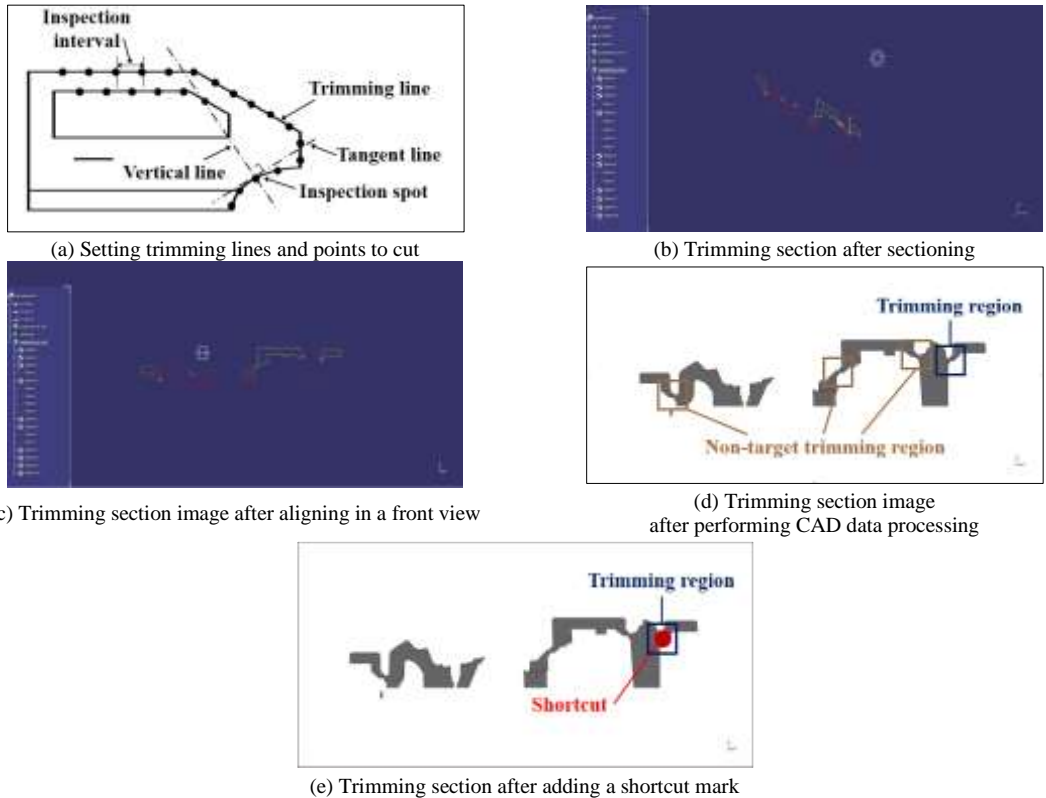

(a) Setting trimming lines and points to cut

(b) Trimming section after sectioning

(c) Trimming section image after aligning in a front view

(d) Trimming section image after performing CAD data processing

(e) Trimming section after adding a shortcut mark

Fig. 6. Intermediate results of data processing

Section images may contain non-target trimming regions along with the target trimming region, as shown in Fig. 6(d). The target trimming region is the one to be inspected. Non-target trimming regions appear because the sectioning line meets several other trimming lines, as shown in Fig. 6(a). To help a CNN model easily distinguish the target trimming region from non-target trimming regions, we add a shortcut, a red circle mark, on the target trimming region when generating section images, as shown in Fig. 6(e). A shortcut is a piece of additional information that helps the AI models easily learn the given tasks (Geirhos et al., 2020; Minderer et al., 2020; Doersch et al., 2015; Negnevitsky, 2005). Because the shape of the shortcut mark is quite different from that of the trimming regions, it is easy for a CNN model to detect it in an image. By adding a shortcut mark at the target trimming region in the section images, our CNN model can easily detect the target trimming region with a high accuracy. The details are presented in Section 4.3.

*3.3. Trimming region detection*

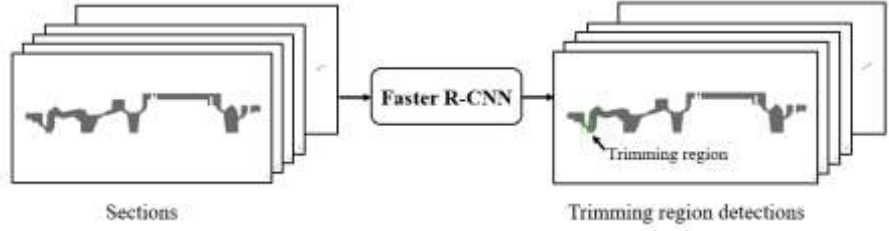

Fig. 7. The process of trimming region detection. Trimming region that output of the target point detection model is too small for the image size.

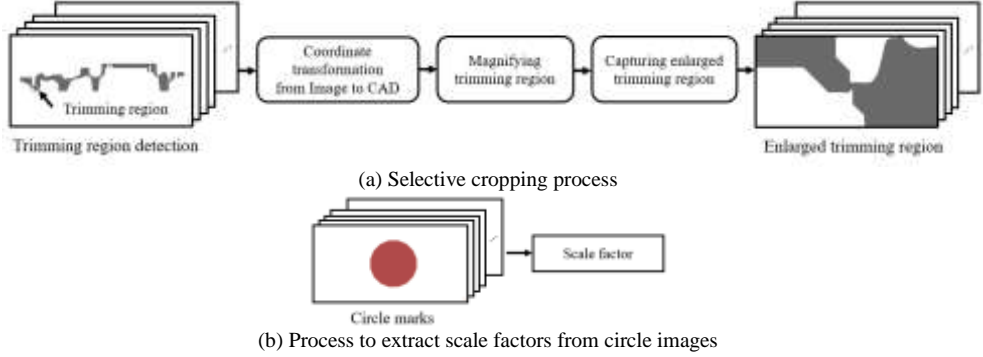

(a) Selective cropping process

(b) Process to extract scale factors from circle images

Fig. 8. The process of selective cropping

In the trimming region detection stage, the target trimming region in the section image is identified to measure the target line lengths, as shown in Fig. 7. We used the Faster R-CNN model for trimming region detection. The input of the model is section images, and the output is the positions of the target trimming region. The engineers generate the training data by manually identifying the positions of the target trimming regions in the section images.

*3.4. Selective cropping*

The size of the target trimming region in the section images is too small: 30 × 30 pixels; thus, it is difficult to measure the target line lengths accurately. We obtain an enlarged image of the target trimming region for a precise measurement in this stage. As shown in Fig. 8(a), the selective cropper generates an enlarged trimming region image using the CAD's built-in magnifying and capturing tools. To obtain enlarged images, we first need to convert the image coordinates of the trimming region into monitor coordinates to perform CAD operations. The coordinates for the trimming region are calculated using a linear relationship to map the CAD interface on a monitor from an image.

The target lengths should be measured in millimeters; however, we can measure the lengths in the section images in pixels. To provide a scale factor on what millimeters a pixel is, we create an additional 20-mm circle image, as shown in Fig. 8(b). The scale factor, $\gamma_{scale}$, is obtained by Eq. 5, where $px_{circle}$ is the measured diameter in pixels, and $d_{circle}$ is the measured diameter in millimeters, as shown in Fig. 9. This scale factor is used in the length measurement stage.

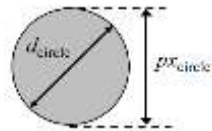

Fig. 9. Circle image. $d_{circle}$ denotes the diameter of a circle (20mm) and $px_{circle}$ denotes the measured pixels of the diameter.

$$\gamma_{scale} = \frac{d_{circle}}{px_{circle}} \qquad (5)$$

*3.5. Target point detection and length measurement*

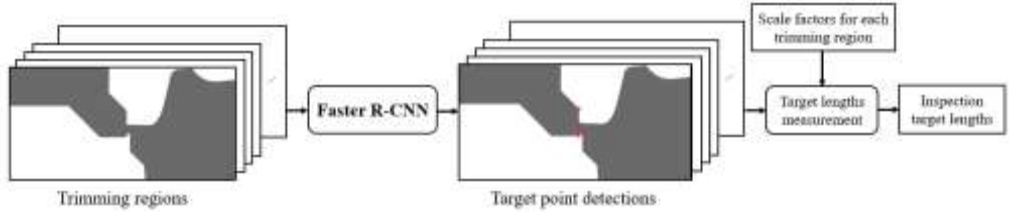

Fig. 10. Intermediate results of target point detection and length measurement

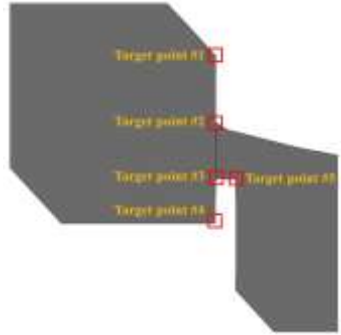

Fig. 11. Trimming region image from selective cropping. The length, measured in pixel, of PD, UL, GL is determined by target points #1~5.

The locations of both endpoints for the target lines in the enlarged trimming region image is identified, as shown in Fig. 10. We use the Faster R-CNN model for target point detection. The input is an enlarged trimming region image, and the outputs are the bounding boxes of the five endpoints, as shown in Fig. 11. The model is trained using the bound boxes of the target endpoints marked by the engineers, as shown in Fig. 11.

We assume a target point is located at the center of the bounding box and measure the distances between the detected target points. The length of the PD is determined by target points #1 and #4, UL by targets #2 and #4, and GL by targets #3 and #5. Because the lengths are measured in pixels, we need to convert them to millimeters using the following equation:

$$l_{tar} = l_{pixel} \times \gamma_{scale} \qquad (6)$$

where $l_{tar}$ denotes the length in millimeters and $l_{pixel}$ denotes the length in pixels. The measurement accuracy is discussed in Section 4.2.

The operations in the AI domain—trimming region detection and length measurement, and the operations in the CAD domain—data processing, and selective cropping, are seamlessly integrated so that the measuring process is a one-click operation without any expert intervention.

**4. Experimental Results and Analysis**

In this section, we present the experimental results to evaluate the performance of the trimming region detection and target length measurement models. We perform a four-fold cross-validation; therefore, the training dataset includes the images generated from three CAD files, and the validation dataset includes the images generated from the other.

Table 1. 4-fold cross validation accuracy on the trimming region detection

| Test Design | Trimming region detection accuracy |
|---|---|
| Design 1 | 100% |
| Design 2 | 100% |
| Design 3 | 100% |
| Design 4 | 100% |
| Average | 100% |

*4.1. Trimming region detection*

Our model accurately detects all trimming regions, even though we have trained it with a small CAD dataset, as shown in Table 1. Because we generate simplified section images by background removal, color reduction, and tree information removal, we can train the CNN model with a small dataset from only three CAD files. Additionally, we add circle marks as shortcuts to the section images, which provides additional information to detect the target trimming regions during training. This prevents the model from confusing the trimming regions with other non-target regions; thus, it can easily detect the trimming regions and be rapidly trained. CAD data processing and shortcut addition are effective for trimming region detection. More details regarding the effect of shortcut addition are discussed in Section 4.4.

*4.2. Target point detection and length measurement*

Table 2. 4-fold cross validation accuracy on target point detection

| Test Design | Target point detection accuracy |
|---|---|
| Design 1 | 97.6% |
| Design 2 | 98.4% |
| Design 3 | 100.0% |
| Design 4 | 97.2% |
| Average | 98.3% |

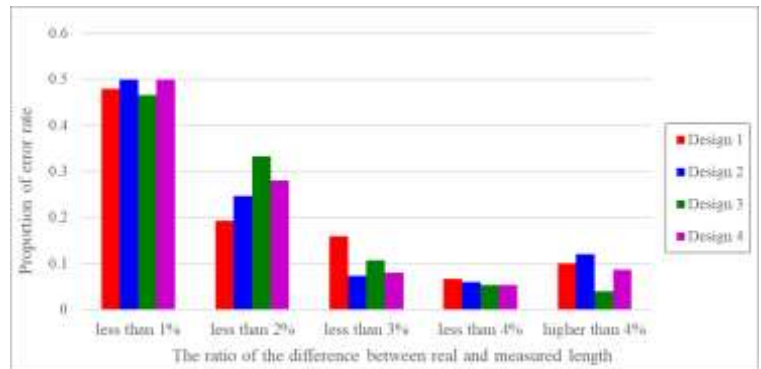

Fig. 12. The proportion of the length difference between real and measured.

The target line endpoints are accurately detected by our system, as shown in Table 2. The error rate of the target point detection is only 1.7 % on an average, which means the system may miss only 1.7 % of the endpoints. Because the trimming regions are identified and enlarged in the selective cropping stage, the input images, which are enlarged around the trimming region, include the characteristics of the endpoints. This allows the target point detection model to effectively train and easily detect the endpoints.

Because the CAD files have different trimming die designs, there are minor differences in the detection accuracy for each validation file. In some cases, the detection accuracy is not 100 %, which means some target points are missed. The probability of detecting all five target points in a trimming region is approximately 91.8% because the detection accuracy is 98.3 % on an average. Thus, the system may fail to process approximately 8.2 % of the section images. We can cope with dense inspections. If the number of inspection spots is doubled, the inspection failure rate would decrease to 0.67 % (8.2 % × 8.2 %). In this program, we have used 5 – 50 spots for inspection of one trimming line, which reduces the failure rate to at least 0.00037% or less. This multi-spot inspection is possible because the program is fully automated. Because an automatic inspection is faster than a manual inspection, we can handle more inspection spots with less time and cost. More details are provided in Section 4.5.

The target lengths, PD, UL, and GL are based on the detected target points and scale factors. Table 3 shows the average value of the length measurement error and Fig. 12 shows its distribution. The length measurement error is the difference between the real and measured lengths, divided by the real length. It is only 2.4 % on an average, and approximately half of the measured target lengths have an error lower than 1 %, as shown in Fig. 12. Approximately 90 % of the measurements have an error lower than 4 %. An average length measurement error of 2.4 % is at the 14th place out of 18 grades, which is suitable for manufacturing operations such as pressing, rolling, and other forming operations according to ISO 286-1:1988 standard tolerance grades (ISO 286-1, 1988). It also passes the quality management team's evaluation in the automotive company.

Table 3. 4-fold cross validation accuracy on target point detection

| Test Design | Average error of length measurement |
|---|---|
| Design 1 | 2.0% |
| Design 2 | 2.7% |
| Design 3 | 1.3% |
| Design 4 | 3.5% |
| Average | 2.4% |

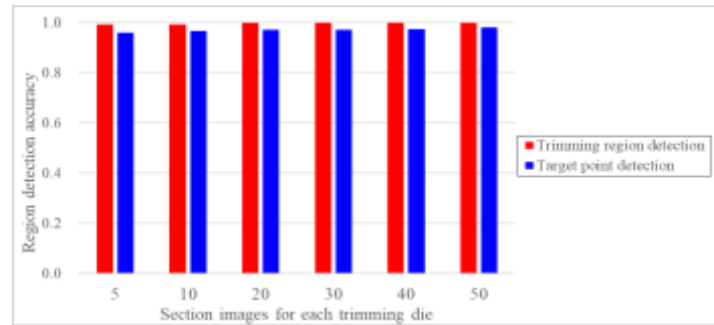
(a) Trimming region and target point detection accuracy

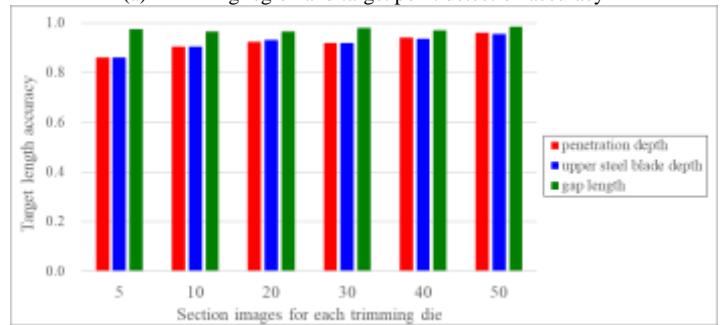
(b) Target length accuracy

Fig. 13. The accuracy of detection model and target length via number of training data.

*4.3. Performance according to the number of section images*

Deep learning techniques are affected by the size of training dataset. Therefore, we verify the sensitivity of our system to the size of training dataset. We observe the accuracy of our model according to the number of section images from each trimming die, as shown in Fig. 13. We generate 5–50 section images for each trimming die design. When five section images are generated, we use 15 section images to train our model; when 50 section images are generated, we use 150 images for training.

Fig. 13(a) shows the trimming region detection and target endpoint detection performance according to the number of images. The region detection accuracies are 99.5 % in the cases of 5 and 10 section images per die and 100 % in the other cases. Although we have trained the model with a small number of section images, the detection performance is quite high. The target point detection results also show the performance improvement with an increasing number of training images, reaching 98.3 % when the model is trained with 150 images. The results show that our proposed data processing, shortcuts, and selective cropping methods can effectively train models with only a small amount of data.

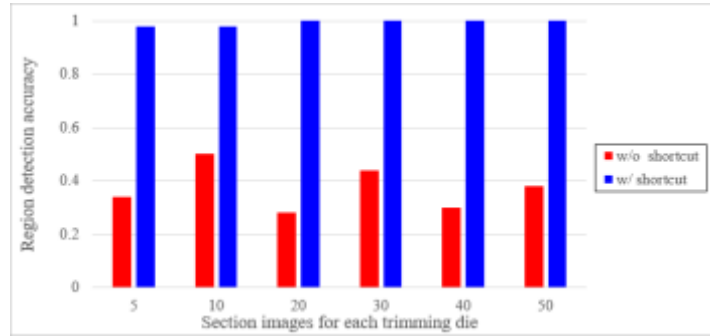

Fig. 14. Trimming region detection accuracy via presence or absence of shortcuts.

Fig. 13(b) shows the accuracy of target length measurement. Length measurement is more complex than the detection of regions and endpoints. We need more training images for a high performance. With five section images per die, the accuracies of PD, UD, and GL are 86.0 %, 86.0 %, and 97.5 %, respectively; however. with 50 section images per die, they are 96.0 % for PD, 95.5 % for UD, and 98.5 % for GL. The performance improves with an increase in the number of images, and the performance with 50 section images from each die is satisfactory.

*4.4. Effect of the existence of shortcuts*

To verify the effect of shortcuts, we observe the trimming region detection accuracy with and without shortcuts, as shown in Fig. 14. Because non-target regions are similar to target trimming regions, it is difficult to distinguish one from the other without shortcuts. The accuracies are lower than 50 % without shortcuts. Even with 50 section images for each trimming die, the accuracy is lower than 40 %. When the model is trained from section images with shortcuts, the detection accuracy is 100 % in most cases. The results show that our shortcut addition method makes it easier for AI models to detect trimming regions, and is an effective way to significantly improve performance.

*4.5. Comparison of inspection time*

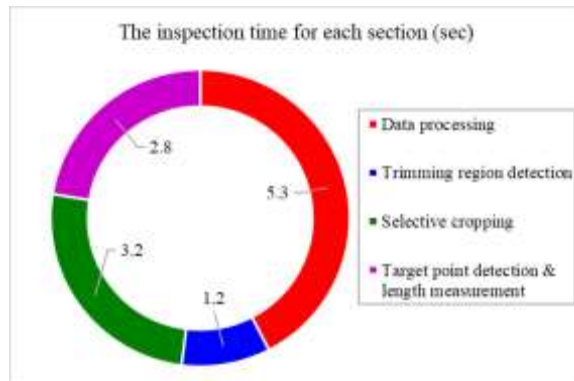

Fig. 15. Inspection time for each spot. It takes 5.3 seconds for data processing, 3.2 seconds for selective cropping, 1.2 seconds for trimming region detection, and 2.8 seconds for point target detection & length measurement.

We also analyze the inspection time of the proposed system. We use a system with an AMD Ryzen 5 3600 6-Core processor and NVIDIA GTX 1660 Ti. We measure the processing time in the AI domain, processing time in the CAD domain, and zigzag interaction time. The AI domain time includes the time required for trimming region detection, target point detection, and length measurement. The CAD domain time includes the time required for data processing and selective cropping, and the zigzag interaction time is the time required to transfer the output between the CAD and AI domains.

The times shown in Fig. 15 are the averages of over 50 inspections. One inspection takes 5.3 seconds for data processing, 3.2 seconds for selective cropping, 1.2 seconds for trimming region detection, and 2.8 seconds for target point detection and target length measurement. The zigzag interaction requires 10 seconds for a trimming die. It is almost independent of the number of inspection spots because we have processed all section images in a batch mode. If there are 50 inspection spots in a trimming die, the total time will be $((5.3+3.2+1.2+2.8) \times 50) +10$ seconds.

The inspection time in the CAD domain is approximately 8.5 seconds, and that in the AI domain is approximately 4.0 seconds. Most inspection times are in the CAD domain. It is 2.1 times longer than that in the AI domain because it involves CAD operations such as capturing and sectioning.

In contrast, engineers take about 60 seconds to inspect three target lengths: PD, UL, and GL. Our automatic inspection program takes approximately 10.5 seconds to inspect one spot. The proposed system reduces the inspection time by 82.5 % when compared with manual inspection. Thus, our system can significantly reduce the time and human resources needed for inspection.

*4.6. Extensibility to various CAD software*

We integrate the AI modules with CATIA, a popular CAD software. Because there are several CAD software, the extensibility of the inspection system is crucial. The system needs to be easily integrated with other CAD software.

In our automatic inspection system, 3D CAD files are converted to 2D images for interaction between the CAD software and AI modules, and the output of the AI modules is transferred to the CAD software through VBA macros, as shown in Fig. 16. It is designed for maximum extensibility to other CAD software. Our system uses functions related to capturing, enlargement, and macro-development with VBA. These functions are available not only in the CATIA program but also in other CAD software, and the VBA used to automate a process is also available in other CAD software. Image interfacing is effective for extending the proposed system to other CAD software.

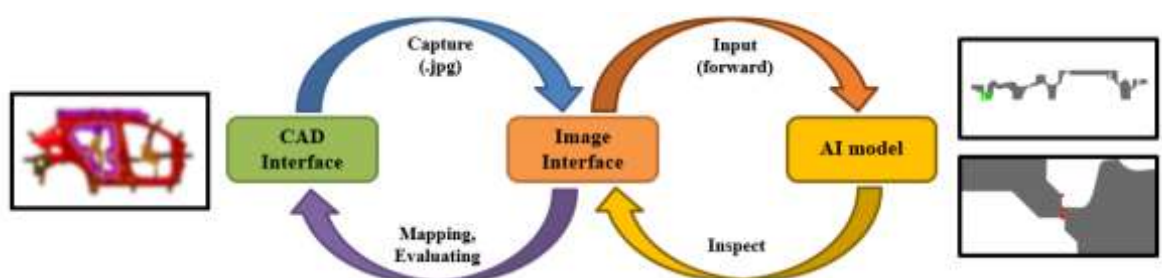

Fig. 16. Interactive our inspection system between CAD Interface and AI model. The trimming die image is partially blurred by the manufacturer's requests.

5. Conclusion

We automated an automobile trimming die design inspection, which is labor-intensive, time-consuming, and prone to human error, by integrating AI modules with CAD software. The AI modules detect the trimming regions and target endpoints, and the CAD software performs the operations requested by the AI modules. Our system is CAD-independent, data-efficient, and easy to use. We designed a zigzag process for a seamless interface and maximum independence between CAD software and AI modules. The AI modules were trained with a few trimming dies by pre-processing 3D CAD files and adding shortcuts. The system has a high usability so that an inspection could be performed with only a single click. Its high inspection performance passed the evaluation of the automotive company's quality management team, and the inspection time was reduced to approximately one-fifth of the manual inspection time by experts. The proposed zigzag process-based design inspection system can be extended to other design inspection tasks using other CAD software.

**CRediT authorship contribution statement**

**Jinsub Lee:** Methodology, Draft, Investigation, Visualization, Data Curation, Validation. **Tae-Hyun Kim:** Methodology, Draft, Investigation, Visualization, Data Curation. **Sang-Hwan Jeon:** Resources, Project administration. **Sung-Hyun Park:** Data Curation, Validation, Resources, Software. **Sang-Hi Kim:** Resources, Project administration. **Eun-Ho Lee:** Draft, Supervision, Conceptualization, Review & Editing. **Jee-Hyong Lee:** Draft, Supervision, Conceptualization, Review & Editing.

**Acknowledgments**

This work was partially supported by Stamping Tool Design Team and Advanced Manufacturing CAE Team at Hyundai motor group, the National Research Foundation of Korea(NRF) grant, funded by the Korea government(MSIT)(No. 2019-0-00421), and the HPC Support Project, supported by the 'Ministry of Science and ICT' and NIPA.

**Declaration of competing interest**

The authors declare that they have no known competing financial interests or personal relationships that could have appeared to influence the work reported in this paper.